\documentclass[11pt]{article}
\newcommand\bi[2]{{{#1}\atopwithdelims(){#2}}}
\newcommand\qbi[3]{{{#1}\atopwithdelims[]{#2}}_{#3}}

\usepackage{epsfig} \usepackage{amssymb,amsfonts,latexsym}
\newtheorem{theo}{Theorem}
 
 \newtheorem{lem}{Lemma}
\newenvironment{proof}{{\em Proof.}}{\mbox{}\hfill $\Box$\medskip}

\newenvironment{rem}{\noindent{\bf Remark. }}{\smallskip}
\newcommand{\la}{\lambda} \title{New Identities of Hall-Littlewood
Polynomials and  Rogers-Ramanujan Type\footnote{Version du 02
oct. 2001} } \author{Fr\'ed\'eric Jouhet and  Jiang
Zeng\\ \small Institut Girard Desargues, Universit\'e Claude
Bernard (Lyon 1)\\ \small 43, bd du 11 Novembre 1918, 69622
Villeurbanne Cedex, France\\ \small E-mail~:
\texttt{\{jouhet,zeng\}@desargues.univ-lyon1.fr}\\ } \date{}
\begin{document}
\maketitle
\begin{abstract}
We prove two new  summation formulae of Hall-Littlewood
polynomials over partitions into bounded parts and derive some new
 multiple $q$-identities  of Rogers-Ramanujan type.
\end{abstract}
\section{Introduction}
The  Rogers-Ramanujan identities~(see \cite{An76,An86})~: $$
\sum_{n=0}^\infty{q^{n^2+an}\over
(1-q)(1-q^2)\cdots(1-q^n)}=\prod_{n=1\atop n\equiv \pm
(a+1)\pmod{5}}^\infty (1-q^n)^{-1}, $$ where $a=0$ or 1, are among
the most famous $q$-series identities in partitions and
combinatorics.
Since their discovery the Rogers-Ramanujan identities have been
proved and generalized in various ways (see \cite{An76,An86, BIS,
St} and the references cited there).
In~\cite{St},  by adapting a
method of Macdonald for calculating partial fraction expansions of
symmetric formal power series, Stembridge gave an unusual proof of
Rogers-Ramanujan identities as well as fourteen other non trivial
$q$-series identities of Rogers-Ramanujan type and their multiple
analogs. Although it is possible to describe his proof  within the
setting of $q$-series, two  summation formulas  of Hall-Littlewood
polynomials were  a crucial source of inspiration for such kind of
identities. One of our original motivations was to look for new multiple
$q$-identities of Rogers-Ramanujan type through this approach, but
we think that the new summation
formulae of Hall-Littlewood polynomials are interesting for their own.

Throughout this paper we will use the standard  notations of
$q$-series (see, for example, \cite{GR}). Set $(x)_0:=(x;q)_0=1$ and
for $n\geq 1$
\begin{eqnarray*}
(x)_n&:=&(x;q)_n=\prod_{k=1}^n(1-xq^{k-1}),\\ (x)_\infty&:=&
(x;q)_\infty=\prod_{k=1}^\infty(1-xq^{k-1}).
\end{eqnarray*}
For $n\geq 0$ and $r\geq 1$, set $$
(a_1,\cdots,a_r;\,q)_n=\prod_{i=1}^r(a_i;q)_n\,\,,\qquad
(a_1,\cdots ,a_r;q)_\infty=\prod_{i=1}^r(a_i;\,q)_\infty. $$

Let $n\geq 1$ be a fixed integer and  $S_n$  the group of
permutations of the set $\{1,\, 2, \ldots, n\}$. Let
$X=\{x_1,\ldots, x_n\}$ be a  set of indeterminates and $q$ a
parameter. For each \emph{partition} $\lambda=(\la_1,\ldots,
\la_n)$ of length $\leq n$, if 
$m_i:=m_i(\la)$ is  the multiplicity
of $i$  in $\la$, then we also note $\la$ by $(1^{m_1}\,2^{m_2}\, \ldots)$. Recall that the Hall-Littlewood polynomials
$P_{\la}(X,q)$ are defined by~\cite[p.208]{Ma}~: 
$$ P_{\la}(X,q)=
\prod_{i\geq 1}\frac{(1-q)^{m_i}}{(q)_{m_i}}\,\sum_{w\in
S_n}w\left(x_1^{\la_1}\ldots x_n^{\la_n}
\prod_{i<j}\frac{x_i-qx_j}{ x_i-x_j}\right), $$ where the factor
is added  to ensure the coefficient of $x_1^{\la_1}\ldots x_n^{\la_n}$ in $P_\la$ is 1.

For a parameter $\alpha$  define  the auxiliary function $$
\Psi_q(X;\alpha):=\prod_{i}(1- x_i)^{-1} (1-\alpha
x_i)^{-1}\prod_{j<k}\frac{1-qx_jx_k}{1-x_jx_k}. $$ Then it is
well-known~\cite[p. 230]{Ma}
 that the sums of $P_{\la}(X,q)$
over all partitions and even partitions are given by the following
formulae~:
\begin{eqnarray}
\sum_{\la}P_\la(X,q)&=&\Psi_q(X;0),\label{lit1}\\
\sum_{\la}P_{2\la}(X,q)&=&\Psi_q(X;-1).\label{lit2}
\end{eqnarray}
For any sequence $\xi\in \{\pm 1\}^n$  set
$X^{\xi}=\{x_1^{\xi_1},\cdots,x_n^{\xi_n}\}$. Then, by summing $P_\la$
over partitions with bounded parts, Macdonald
\cite[p. 232]{Ma} and Stembridge~\cite{St} have respectively
generalized (\ref{lit1}) and (\ref{lit2}) as follows~:
\begin{eqnarray}
\sum_{\la_1\leq k}P_\la(X,q)&=&\sum_{\xi\in \{\pm
1\}^n}\Psi_q(X^{\xi};0)\prod_i x_i^{k(1-\xi_i)/2},\label{macdo}\\
\sum_{\la_1\leq 2k \atop \la\,\, even}P_{\la}(X,q)&=&\sum_{\xi\in
\{\pm 1\}^n}\Psi_q(X^{\xi};-1)\prod_i
x_i^{k(1-\xi_i)}.\label{stem}
\end{eqnarray}
Now, for parameters $\alpha$, $\beta$
define another auxiliary function 
$$
\Phi_q(X;\alpha,\beta):=\prod_{i} \frac{1- \alpha x_i}{1-\beta
x_i}\,\prod_{j<k}\frac{1-qx_jx_k}{1-x_jx_k}.
 $$ 
Then the following summation formulae similar to (1) and (2) for
Hall-Littlewood polynomials  hold true~\cite[p.232]{Ma}~:
\begin{eqnarray}
\sum_{\la'\,\textrm{even}}c_{\la}(q)\,P_\la(X,q)&=&\Phi_q(X;0,0)\label{f},\\
\sum_{\la}d_{\la}(q)\,P_\la(X,q)&=&\Phi_q(X;q,1)\label{g},
\end{eqnarray}
 where $\la'$ is the conjugate of $\la$ and
$$ c_{\la}(q)=\prod_{i\geq 1}(q;q^2)_{m_i(\la)/2},\qquad
d_\la(q)=\prod_{i\geq 1}\frac{(q)_{m_i(\la)}}{(q^2;q^2)_{[m_i(\la)/2]}}. $$

In view of the numerous applications of (3) and (4) it is natural to seek such
extensions for  (5) and (6).
 However,
as remarked by Stembridge~\cite[p. 475]{St}, in these other cases
there arise complications which render \emph{doubtful} the
existence of expansions as explicit as those of (\ref{macdo}) and
(\ref{stem}). We noticed that these complications arise if one
wants to keep exactly the same coefficients $c_\la(q)$ and $d_\la(q)$
as in (5) and (6) for the sums over bounded partitions.
Actually we have the following
\begin{theo} For $k\geq 1$,
\begin{eqnarray}
 \sum_{\la_1\leq k \atop \la'\,\,
even}c_{\la,k}(q)P_{\la}(X,q)&=&\sum_{\xi\in \{\pm 1\}^n\atop
|\xi|_{-1}\,\, even}\Phi_q(X^{\xi};0,0)\prod_i x_i^{k(1-\xi_i)/2},
\label{moi1}\\
 \sum_{\la_1\leq k}d_{\la,k}(q)P_{\la}(X,q)&=&\sum_{\xi\in \{\pm
1\}^n}\Phi_q(X^{\xi};q, 1)\prod_i
x_i^{k(1-\xi_i)/2},\label{autre}
\end{eqnarray}
where $$ c_{\la,k}(q)=\prod_{i=1}^{k-1}(q;q^2)_{m_i(\la)},\quad
d_{\la,k}(q)=\prod_{i=1}^{k-1}
\frac{(q)_{m_i(\la)}}{(q^2;q^2)_{[m_i(\la)/2]}}.
 $$
\end{theo}

\begin{rem}  We were led to such extensions by
starting  from the right-hand side instead of the left-hand side
and inspired by the  similar formulae corresponding to the case
$q=0$  of Hall-Littlewood polynomials~\cite{JZ}, i.e., Schur
functions. In the initial stage we made also the Maple tests using
the package ACE~\cite{ACE}. In the  case $q=0$,  the right-hand
sides of (3), (4), (7)  and (8) can be written as quotients of
determinants and the formulae reduce to the known identities of
Schur functions~\cite{JZ}.
\end{rem}

For any partition $\la$ it will be convenient  to adopt  the
following notation~: $$(x)_\la:= (x;\,
q)_\la=(x)_{\la_1-\la_2}(x)_{\la_2-\la_3}\cdots, $$ and to
introduce the general $q$-binomial coefficients $$
\qbi{n}{\la}{}:=\frac{(q)_n}{(q)_{n-\la_1}(q)_{\la}}, $$ with the
convention that $\qbi{n}{\la}{}=0$ if $\la_1>n$. If $\la=(\la_1)$
we recover the classical $q$-binomial coefficient. Finally, for
any partition $\la$ we denote by $l(\la)$ the length of $\la$,
i.e., the number of its positive parts, and
$n(\la):=\sum_i\left({\la_i\atop 2}\right)$.

The  following  is the key $q$-identity which allows to produce
identities of Rogers-Ramanujan type.
\begin{theo} For $k\geq 1$,
\begin{eqnarray}
&&\hskip -0.5cm \sum_{l(\la)\leq k}z^{|\la|}q^{n(2\la)}\frac{(a,\,
b;q^{-2})_{\la_1}}{(q^2;q^2)_\la(q;q^2)_{\la_k}}=\frac{(z;q^2)_\infty}{(abzq
;q^2)_\infty}\label{moi5}\\ &&\times \sum_{r\geq
0}z^{kr}q^{(k+1)\left({2r \atop
2}\right)}\frac{(a,b;q^{-2})_r(aq^{2r+1}z,bq^{2r+1}z;q^2)_\infty}{(q)_{2r}
(zq^{2r-1})_\infty}(1-zq^{4r-1})\nonumber.
\end{eqnarray}
\end{theo}

The remainder of this paper is organized as follows~:  in section~2 we first
give multiple analogs of Rogers-Ramanujan type identities which
are consequences of Theorem~2, in section 3 we give the proof of
Theorem 1 and some consequences, and defer the
elementary proof, i.e., without using the Hall-Littlewood polynomials,
of Theorem~2 and other multiple $q$-series
identities to section 4.
 In section~5 we will compare our multianalogs of Rogers-Ramanujan's
type identities with those obtained through Andrews-Bailey's method.

\section{Multiple identities of Rogers-Ramanujan type}
We need the {\sl Jacobi triple product}
identity~\cite[p.21]{An76}~:
\begin{equation}
J(x,\, q):=1+\sum_{r= 1}^\infty (-1)^rx^rq^{\left({r \atop
2}\right)}(1+q^r/x^{2r})=(q,\, x,\,  q/x;\, q)_\infty. \label{Ja2}
\end{equation}
For any partition $\la$ set
 $n_2(\la)=\sum_i \la_i^2$.
We derive then from Theorem~2 the following 
identities  of Rogers-Ramanujan type.
\begin{theo} For $k\geq 1$,
\begin{equation}\label{moi2}
\sum_{l(\la)\leq k}
 \frac{q^{2n_2(\lambda) }}
{(q;q^2)_{\lambda_k}(q^2;q^2)_{\lambda}} =\prod_n (1-q^n)^{-1}
\end{equation}
where $n\equiv \pm(2k+1),\, \pm(2k+3),\pm2, \pm4,\,\ldots, \pm4k
\pmod{8k+8}$;
\begin{equation}\label{moi6}
\sum_{l(\la)\leq k}
 \frac{q^{2n_2(\lambda)-2\la_1}}
{(q;q^2)_{\lambda_k}(q^2;q^2)_{\lambda}}(1-q^{2\la_1})
=\frac{(q^{2k-1},\, q^{6k+9};\;q^{8k+8})_\infty}{\prod_n (1-q^n)}
\end{equation}
where $n\equiv  \pm(2k+5), \pm2,\,\ldots, \pm4k, \pm(4k+2)
\pmod{8k+8}$;
\begin{eqnarray}
&&\hskip -1 cm\sum_{l(\la)\leq k}\frac{q^{2n_2(\lambda)-\la_1^2}}
{(q;q^2)_{\lambda_k}(q^2;q^2)_{\lambda}}(-q;q^2)_{\la_1}
\nonumber\\ &&\hskip 2cm=\frac{(-q;q^2)_\infty}{(q^2;q^2)_\infty}
(q^{4k+2},\, -q^{2k},-q^{2k+2};\;q^{4k+2})_\infty\label{moi7};\\
&&\hskip -1 cm\sum_{l(\la)\leq k}\frac{q^{2n_2(\lambda)-
\la_1^2-\la_1}}{(q;q^2)_{\lambda_k}(q^2;q^2)_{\lambda}}(-1;q^2
)_{\la_1}(1-q^{2\la_1})\nonumber\\ &&\hskip
2cm=\frac{(-q^2;q^2)_\infty}{(q^2;q^2)_\infty}(q^{4k+2},\,-q^{
2k-1},\,-q^{2k+3};\;q^{4k+2})_\infty;\label{moi8}\\ &&\hskip -1
cm\sum_{l(\la)\leq k}\frac{q^{2n_2(\lambda)
-\lambda_1^2+\la_1}}{(q;q^2)_{\lambda_k}(q^2;q^2)_{\lambda}}(-1;q^2)_{\la_1
}(-q;q^2)_{\la_1}\nonumber\\ &&\hskip
2cm=\frac{(-q)_\infty}{(q)_\infty} (q^{4k}, -q^{2k},
-q^{2k};q^{4k})_\infty;\label{moi12}\\ &&\hskip -1
cm\sum_{l(\la)\leq
k}\frac{q^{2n_2(\lambda)-\la_1^2+\la_1}}{(q;q^2)_{\la_k}(q^2;q^2)_{\la}}(-1;q^2)_{\la_1}
\nonumber\\ &&\hskip 2cm=
\frac{(-q^2;q^2)_\infty}{(q^2;q^2)_\infty} (q^{4k+2},\,
-q^{2k+1},\, -q^{2k+1};\; q^{4k+2})_\infty.\label{moi11}
\end{eqnarray}
\end{theo}
\begin{proof}
Set $z=q$ in (\ref{moi5}),
\begin{eqnarray}
&&\sum_{l(\la)\leq k}q^{|\la|+n(2\la)}\frac{(a,\, b;\,
q^{-2})_{\la_1}}{(q^2;\, q^2)_\la(q;\, q^2)_{\la_k}}=\label{zq}\\
&&\hskip 1cm 1+\sum_{r\geq 1}q^{2kr^2+\bi{2r}{2}} {(a,b;\,
q^{-2})_r(aq^{2r+2}, bq^{2r+2};\, q^2)_\infty \over (abq^2;\,
q^2)_\infty(q^2;\, q^2)_\infty}(1+q^{2r}).\nonumber
\end{eqnarray}

For (11), setting  $a=b=0$ in (\ref{zq}) we obtain
$$ \sum_{l(\la)\leq k}
 \frac{q^{2n_2(\lambda) }}
{(q;q^2)_{\lambda_k}(q^2;q^2)_{\lambda}}
=(q^2;q^2)_\infty^{-1}J(-q^{2k+1},\, q^{4k+4}). 
$$ 
The
right side of (\ref{moi2}) follows then from (\ref{Ja2})  after simple
manipulations.

For (12), set $a=0$ in (\ref{zq}) and  multiply  both sides by
$1-q^{-2}$.  Identifying the  coefficients of $b$ we obtain~:
$$
\sum_{l(\la)\leq k}\frac{q^{2n_2(\lambda)-2\la_1}}
{(q;q^2)_{\lambda_k}(q^2;q^2)_{\lambda}}(1-q^{2\la_1})
=(q^2;q^2)_\infty^{-1}\,J(-q^{2k-1};\,q^{4k+4}).
$$
The result  follows from (\ref{Ja2}) after simple manipulations.

Identity (\ref{moi7}) follows from  (\ref{zq}) with $a=-q^{-1}$
and $b=0$ and then by applying (\ref{Ja2}) with $q$ replaced by
$q^{4k+2}$ and $x=-q^{2k}$.

 For (\ref{moi8}), we choose $a=-1$ in
(\ref{zq}) and multiply both sides by $1-q^{-2}$, then identify
the
 coefficient of $b$. The  identity follows then by applying
 (\ref{Ja2}) with $q$
replaced by $q^{4k+2}$ and $x=-q^{2k-1}$.

Identity (\ref{moi12}) follows from  (\ref{zq})  by taking
$a=-q^{-1}$ and $b=-1$ and then applying (\ref{Ja2}) with $q$
replaced by $q^{4k}$ and $x=-q^{2k}$. For (\ref{moi11}), we choose
$a=-1$ and $b=0$ in (\ref{zq}). The  identity follows then by
applying
 (\ref{Ja2}) with $q$
replaced by $q^{4k+2}$ and $x=-q^{2k+1}$.
\end{proof}

When $k=1$ the above six identities reduce respectively to the
following Rogers-Ramanujan type identities~:
\begin{eqnarray}
\sum_{n=0}^\infty{q^{2n^2}\over (q)_{2n}}&=&\prod_{n=1\atop
n\equiv\pm 2,\pm 3, \pm 4,\pm 5 \pmod{16}}^\infty \frac{1}{1-q^n}
\label{moi3},\\ \sum_{n=0}^\infty{q^{2n^2+2n}\over (q)_{2n+1}}
&=&\prod_{n=1\atop n\equiv\pm 1,\pm 4, \pm 6,\pm 7
\pmod{16}}^\infty \frac{1}{1-q^n} \label{moi4},\\
\sum_{n=0}^\infty q^{n^2}{(-q;q^2)_n\over
(q)_{2n}}&=&{(q^6,q^6,q^{12};\, q^{12})_\infty\over (q)_\infty},
\label{moi9}\\ \sum_{n=0}^\infty q^{n^2+n}{(-q^2;q^2)_n\over
(q)_{2n+1}}&=&{(q^3,q^9,q^{12};\, q^{12})_\infty\over (q)_\infty},
\label{moi10}\\ 
1+2\sum_{n\geq
1}q^n\frac{(-q)_{2n-1}}{(q)_{2n}}&=&\frac{
(q^{4},\,-q^{2},\,  -q^{2};\,q^{4})_\infty}{(q)_\infty(q;q^2)_\infty}\label{moi14},\\
1+2\sum_{n\geq
1}q^{n(n+1)}\frac{(-q^2;q^2)_{n-1}}{(q)_{2n}}&=&
\frac{ (q^{6}, -q^{3}, -q^{3};\,
q^{6})_\infty}{(q)_\infty\,(-q;q^2)_\infty}.\label{moi13}
\end{eqnarray}

Note that (\ref{moi3}), (\ref{moi4}), (\ref{moi9}) and
(\ref{moi10}) are already known, they correspond to Eqs.~(39), (38), (29) and (28) in Slater's  list \cite{Sl2},
respectively, but (\ref{moi14}) and (\ref{moi13}) seem  to be new.
\section{Proof of Theorem~1 and consequences}
\subsection{Proof of identity (\ref{moi1})}
For any statement $A$ it will be convenient to use the true or false function
 $\chi(A)$, which is 1 if $A$ is true
 and 0 if $A$ is false.
Consider the generating function $$
S(u)=\sum_{\la_0,\la}\chi(\la'\,\textrm{even})\,c_{\la,\la_0}(q)P_\la(X,q)\,
u^{\la_0} $$ where the sum 
is over all partitions $\la=(\la_1,\ldots, \la_n)$ and the integers
 $\la_0\geq \la_1$. Suppose
$\la=(\mu_1^{r_1}\, \mu_2^{r_2}\, \ldots
\mu_k^{r_k})$,
 where 
$\mu_1>\mu_2>\cdots >\mu_k\geq 0$ and $(r_1, \ldots, r_k)$ is a
 composition of $n$.

Let $S_n^\la$ be the set of permutations of $S_n$ which fix $\la$.
 Each $w\in S_n/S_n^\la$ corresponds to  a surjective mapping $f:
X\longrightarrow \{1,2,\ldots, k\}$ such that $|f^{-1}(i)|=r_i$.
For any subset $Y$ of $X$, let $p(Y)$ denote  the product of the
elements of $Y$ (in particular, $p(\emptyset)=1$). We can rewrite
Hall-Littlewood functions as follows : $$ P_\la(X,q)=\sum_{f}
p(f^{-1}(1))^{\mu_1}\cdots p(f^{-1}(k))^{\mu_k}
\prod_{f(x_i)<f(x_j)}{x_i-qx_j\over x_i-x_j}, $$ summed over all
surjective mappings $f: X\longrightarrow \{1,2,\ldots, k\}$ such
that $|f^{-1}(i)|=r_i$. Furthermore, each such $f$ determines a
\emph{filtration}  of $X$ : 
\begin{equation}\label{filtration}
 {\cal F}: \quad
\emptyset=F_0\subsetneq F_1\subsetneq \cdots \subsetneq F_k=X, 
\end{equation}
according to the rule $x_i\in F_l\Longleftrightarrow f(x_i)\leq l$
for $1\leq l\leq k$. Conversely, such a filtration ${\cal
F}=(F_0,\, F_1, \ldots, F_k)$ determines a surjection $f:
X\longrightarrow \{1,2,\ldots, k\}$ uniquely. Thus we can write :
\begin{equation}\label{filtre}
P_\la(X,q)=\sum_{\cal F}\pi_{\cal F}\prod_{1\leq i\leq
k}p(F_i\setminus F_{i-1})^{\mu_i},
\end{equation}
summed  over all the filtrations $\cal F$  such that
 $|F_i|=r_1+r_2+\cdots +r_i$ for $1\leq i\leq k$, and
$$ \pi_{\cal F}=\prod_{f(x_i)<f(x_j)}{x_i-qx_j\over x_i-x_j}, $$
where $f$ is the function defined by $\cal F$.

Now let $\nu_i=\mu_i-\mu_{i+1}$ if $1\leq i\leq k-1$ and
$\nu_k=\mu_k$, thus $\nu_i>0$ if $i<k$ and $\nu_k\geq 0$. Since
the  lengths of columns of $\la$ are
 $|F_j|=r_1+\cdots +r_j$
with multiplicities $\nu_j$ for $1\leq j\leq k$, we have
\begin{equation}\label{mac0}
\chi(\la'\,\textrm{even})=\prod_j\chi(|F_j|\hbox{\scriptsize
even}).
\end{equation}
A filtration ${\cal F}$ is called \emph{even} if $|F_j|$ is even for
$j\geq 1$.
Furthermore, let  $\mu_0=\la_0$ and $\nu_0=\mu_0-\mu_1$ in the
definition of $S(u)$, so that $\nu_0\geq 0$ and
$\mu_0=\nu_0+\nu_1+\cdots +\nu_k$. Define
$\varphi_{2n}(q)=(1-q)(1-q^3)\cdots(1-q^{2n-1})$ and 
 $c_{\cal F}(q)=\prod_{i=1}^k\varphi_{|F_i\setminus
F_{i-1}|}(q)$ for even filtrations $\cal F$. Thus, since $r_j=m_{\mu_j}(\la)$ for $ j\geq 1$,
we have
\begin{eqnarray*}
c_{\la,\,\,\la_0}(q)&=&c_{\cal F}(q) \left(\chi(\nu_k=0)
\varphi_{|F_k\setminus F_{k-1}|}(q)+\chi(\nu_k\neq 0)\right
)^{-1}\\ &&\hspace{2cm}\times \left
(\chi(\nu_0=0)\varphi_{|F_1|}(q)+\chi(\nu_0\neq 0)\right)^{-1}.
\end{eqnarray*}
Let $F(X)$ be the set of filtrations of $X$. Summarizing we obtain 
\begin{eqnarray}
S(u)&=&\sum_{{\cal F}\in F(X)}c_{\cal F}\,\pi_{\cal F}\, \chi({\cal
F\hbox{\scriptsize
even}})\sum_{\nu_1,...,\nu_{k-1}>0}u^{\nu_j}p(F_j)^{\nu_j}\nonumber\\
&\times&\sum_{\nu_0\geq 0}{u^{\nu_0}\over
\chi(\nu_0=0)\,\varphi_{|F_1|}(q) +\chi(\nu_0\neq 0)}\nonumber\\ 
&\times&\sum_{\nu_k\geq 0}
{u^{\nu_k}\,p(F_k)^{\nu_k}\over \chi(\nu_k=0)\,\varphi_{|F_k\setminus
F_{k-1}|}(q)+\chi(\nu_k\neq 0)}\label{bigmac'}.
\end{eqnarray}
For any filtration $\cal F$ of $X$ set 
$$ {\cal A}_{\cal F}(X,
u)=c_{\cal F}(q)\,\prod_{|F_j|\, \hbox{\scriptsize even}}\left[{p(F_j)u \over
1-p(F_j)u}+{\chi(F_j=X)\over \varphi_{|F_j\setminus
F_{j-1}|}(q)}+{\chi(F_j=\emptyset)\over
\varphi_{|F_1|}(q)}\right], 
$$
if ${\cal F}$ is even and 0 otherwise. It follows from (\ref{bigmac'}) that
$$ 
S(u)=\sum_{{\cal F}\in F(X)}\pi_{\cal
F}{\cal A}_{\cal F}(X, u). 
$$ 
Hence  $S(u)$ is a rational function of $u$ with simple
 poles at  $1/p(Y)$, where $Y$ is a subset of $ X$ such that $|Y|$ is
 even. 
 We are
now proceeding to compute  the corresponding  residue $c(Y)$
 at each pole $u=1/p(Y)$.

Let us  start with
$c(\emptyset)$. Writing $\la_0=\la_1+k$ with $k\geq 0$, we see
that
\begin{eqnarray*}
S(u)&=&\sum_\la \chi(\la'\,\textrm{even})\,
c_\la(q)P_\la(X,q)u^{\la_1}\sum_{k\geq 0}{u^k\over
\chi(k=0)\varphi_{m_{\la_1}}(q)+\chi(k\neq 0)}\\
 &=&\sum_\la
\chi(\la'\,\textrm{even})\,
c_\la(q)P_\la(X,q)u^{\la_1}\left({u\over 1-u}+{1\over
\varphi_{m_{\la_1}}(q)}\right).
\end{eqnarray*}
It follows from (\ref{f}) that $$
c(\emptyset)=\left[S(u)(1-u)\right]_{u=1}=\Phi_q(X;0,0).
 $$
For the computations of other  residues, we need some more
notations. For any $Y\subseteq X$, let $Y'=X\setminus Y$ and
$-Y=\{x_i^{-1}:x_i\in Y\}$. Let $Y\subseteq X$ such that $|Y|$ is
even. Then
\begin{equation}\label{residu}
c(Y)=\left[\sum_{\cal F}\pi_{\cal F}{\cal A}_{\cal
F}(X;u)(1-p(Y)u)\right]_{u=p(-Y)}.
\end{equation}
If $Y\notin\cal F$, the corresponding summand is  equal to  0.
Thus we need only to consider the following filtrations ${\cal F}$
: $$ \emptyset=F_0\subsetneq \cdots \subsetneq F_t=Y\subsetneq
\cdots \subsetneq  F_k=X\qquad 1\leq t\leq k. $$ We may then split
$\cal F$ into two filtrations ${\cal F}_1$ and ${\cal F}_2$~ :
\begin{eqnarray*}
{\cal F}_1&:& \emptyset \subsetneq  -(Y\setminus
F_{t-1})\subsetneq \cdots \subsetneq -(Y\setminus F_1)\subsetneq
-Y,\\ {\cal F}_2&:& \emptyset \subsetneq F_{t+1}\setminus
Y\subsetneq   \cdots \subsetneq F_{k-1}\setminus Y\subsetneq  Y'.
\end{eqnarray*}
Then, writing  $v=p(Y)u$ and $c_{\cal F}=c_{{\cal F}_1}\times
c_{{\cal F}_2}$, we have
 $$ \pi_{\cal F}(X)=\pi_{{\cal F}_1}(-Y)\pi_{{\cal
F}_2}(Y')\prod_{x_i\in Y, x_j\in
Y'}\frac{1-qx_i^{-1}x_j}{1-x_i^{-1}x_j}, $$ and ${\cal A}_{\cal
F}(X;u)(1-p(Y)u)$ is equal to
\begin{eqnarray*}
&&{\cal A}_{{\cal F}_1}(-Y;v) {\cal A}_{{\cal
F}_2}(Y';v)(1-v)\left(\frac{v}{1-v}+\frac{\chi(Y=X)}{\varphi_{|Y\setminus
F_{t-1}|}}\right)\\ &&\times
\left({\frac{v}{1-v}+\frac{1}{\varphi_{|Y\setminus
F_{t-1}|}(q)}}\right)^{-1}\left({\frac{v}{1-v}
+\frac{1}{\varphi_{|F _{t+1}\setminus Y|}(q)}}\right)^{-1}.
\end{eqnarray*}
Thus when $u=p(-Y)$, i.e., $v=1$,
\begin{eqnarray*}
&&\hskip -1cm\left[{\cal A}_{\cal
F}(X;u)(1-p(Y)u)\right]_{u=p(-Y)}=\\ &&\left[{\cal A}_{{\cal
F}_1}(-Y;v)(1-v) {\cal A}_{{\cal
F}_2}(Y';v)(1-v)\right]_{v=1}\times\prod_{x_i\in Y, x_j\in
Y'}\frac{1-qx_i^{-1}x_j}{1-x_i^{-1}x_j}.
\end{eqnarray*}
Using (\ref{residu}) and the result of $c(\emptyset)$, which can
be written $$ \left[\sum_{\cal F}\pi_{\cal F}{\cal A}_{\cal F}(X,
u)(1-u)\right]_{u=1}=\Phi_q(X;0,0), $$
 we get
$$ c(Y)=\Phi_q(-Y;0,0)\Phi_q(Y';0,0)\prod_{x_i\in Y, x_j\in
Y'}\frac{1-qx_i^{-1}x_j}{1-x_i^{-1}x_j}. $$ Each subset $Y$ of $X$
can be encoded by  a sequence $\xi\in \{\pm 1\}^n$ according to
the rule~: $\xi_i=1$ if $x_i\notin Y$ and $\xi_i=-1$ if $x_i\in
Y$. Hence $$ c(Y)=\Phi_q(X^{\xi};0,0). 
$$ 
Note also that 
$$
p(Y)=\prod_i x_i^{(1-\xi_i)/2},\qquad p(-Y)=\prod_i
x_i^{(\xi_i-1)/2}. 
$$ 
Now, extracting the
coefficients of $u^k$ in the equation~:
$$ S(u)=\sum_{{Y\subseteq X}\atop |Y| \,
\hbox{\scriptsize even}>0} {c(Y)\over 1-p(Y)u},$$
yields 
$$ \sum_{\la_1\leq k \atop
\la'\,\, even}c_{\la,k}(q)P_{\la}(X,q)=\sum_{{Y\subseteq X}\atop
|Y| \, \hbox{\scriptsize even}}c(Y)p(Y)^k. $$ Finally,
substituting  the value of $c(Y)$ in the above formula we obtain
(\ref{moi1}).

\begin{rem}
Stembridge's formula  (\ref{stem}) can
 be derived from Macdonald's (\ref{macdo}) and
Pieri's formula for Hall-Littlewood polynomials. Indeed, one of
Pieri's formulas states that \cite[p. 215]{Ma} :
\begin{equation}
P_\mu(X,q)e_m(X)=\sum_{\la}\prod_{i\geq
1}\left[{\la'_i-\la'_{i+1}\atop
\la'_i-\mu'_i}\right]P_\la(X,q),\label{p1}
\end{equation}
where the sum is over all  partitions $\la$ such that  
 $\mu\subseteq\la$ with $|\la/\mu|=m$ and there is at
most one cell in each row of the Ferrers diagram of $\la/\mu$. It
follows from (\ref{p1}) that $$\sum_{\mu_1\leq 2k\atop
\mu\,\,even}P_{\mu}(X,q)\sum_{m\geq 0}e_m(X)=\sum_{\la_1\leq
2k+1}P_{\la}(X,q),$$ noticing that $\la$ determines in a unique
way $\mu$ even by deleting a cell in each odd part of $\la$, and
thus $\left[{\la'_i-\la'_{i+1}\atop \la'_i-\mu'_i}\right]=1$.
Finally we obtain the result, using the fact that $\prod_i
(1+x_i^{\xi_i})^{-1}=\prod_i (1+x_i)^{-1}\times\prod_i
x_i^{(1-\xi_i)/2}$. It would be interesting to  give a similar
proof of (7) using (\ref{macdo}) and another Pieri
formula~\cite[p. 218]{Ma}.
\end{rem}

\subsection{Proof of identity (\ref{autre})}
As in the proof of (\ref{moi1}), we compute  the generating
function 
$$ F(u)=\sum_{\la_0,\la}d_{\la,\la_0}(q)P_\la(X;q)\,
u^{\la_0} 
$$ where the sum is over all partitions $\la=(\la_1,\ldots, \la_n)$
and integers $\la_0\geq \la_1$. For any
filtration $ {\cal F}$  of $X$ (cf. (\ref{filtration})) set 
$$ d_{\cal
F}(q)=\prod_{i=1}^k\psi_{|F_i\setminus F_{i-1}|}(q),\quad
\hbox{where}\quad \psi_{n}(q)=(q)_n\,
\prod_{j=1}^{[n/2]}(1-q^{2j})^{-1}.
$$ Thus, as
$r_j=m_{\mu_j}(\la)$, $j\geq 1$, we have
\begin{eqnarray*}
d_{\la,\,\,\la_0}(q)&=&d_{\cal F}(q) \left(\chi(\nu_k=0)
\psi_{|F_k\setminus F_{k-1}|}(q)+\chi(\nu_k\neq 0)\right )^{-1}\\
&&\hspace{2cm}\times \left
(\chi(\nu_0=0)\psi_{|F_1|}(q)+\chi(\nu_0\neq 0)\right)^{-1}.
\end{eqnarray*}
In view of (\ref{filtre}) we have
$$ F(u)=\sum_{{\cal F}\in F(X)}\pi_{\cal F}{\cal B}_{\cal F}(X,
u), 
$$ 
where 
 $$
{\cal B}_{\cal F}(X, u)=d_{\cal F}\prod_j\left[{p(F_j)u \over
1-p(F_j)u}+{\chi(F_j=X)\over \psi_{|F_j\setminus
F_{j-1}|}(q)}+{\chi(F_j=\emptyset)\over \psi_{|F_1|}(q)}\right].
$$ It follows that  $F(u)$ is a rational function of $u$ and can
be written as~: $$ F(u)={c(\emptyset)\over 1-u}+\sum_{{Y\subseteq
X}\atop |Y|>0} {c(Y)\over 1-p(Y)u}. $$ Extracting the coefficient
of $u^k$ in the above identity  yields
\begin{equation}\label{rational}
\sum_{\la_1\leq k}d_{\la,k}(q)P_{\la}(X,q)=\sum_{Y\subseteq
X}c(Y)p(Y)^k.
\end{equation}
It remains to compute the residues. Writing $\la_0=\la_1+r$ with
$r\geq 0$, then
\begin{eqnarray*}
F(u)&=&\sum_\la d_\la(q)P_\la(X,q)u^{\la_1}\sum_{r\geq 0}{u^r\over
\chi(r=0)\psi_{m_{\la_1}}(q)+\chi(r\neq 0)}\\
 &=&\sum_\la
d_\la(q)P_\la(X,q)u^{\la_1}\left({u\over 1-u}+{1\over
\psi_{m_{\la_1}}(q)}\right),
\end{eqnarray*}
it follows from (\ref{g}) that
\begin{equation}\label{star}
c(\emptyset)=\left(F(u)(1-u)\right)|_{u=1}=\Phi_q(X;q,1).
\end{equation}
For computations of the other  residues, set $Y'=X\setminus Y$ and 
define, for $Y=F_t$, the two
filtrations~:
\begin{eqnarray*}
{\cal F}_1&:& \emptyset \subsetneq  -(Y\setminus
F_{t-1})\subsetneq \cdots \subsetneq -(Y\setminus F_1)\subsetneq
-Y,\\ {\cal F}_2&:& \emptyset \subsetneq F_{t+1}\setminus
Y\subsetneq   \cdots \subsetneq F_{k-1}\setminus Y\subsetneq  Y'.
\end{eqnarray*}
Then, writing  $v=p(Y)u$ and $d_{\cal F}=d_{{\cal F}_1}\times
d_{{\cal F}_2}$, we have
 $$ \pi_{\cal F}(X)=\pi_{{\cal F}_1}(-Y)\pi_{{\cal
F}_2}(Y')\prod_{x_i\in Y, x_j\in
Y'}\frac{1-qx_i^{-1}x_j}{1-x_i^{-1}x_j}, $$ and ${\cal B}_{\cal
F}(X;u)(1-p(Y)u)$ can be written as
\begin{eqnarray*}
&&{\cal B}_{{\cal F}_1}(-Y;v) {\cal B}_{{\cal
F}_2}(Y';v)(1-v)\left(\frac{v}{1-v}+\frac{\chi(Y=X)}{\psi_{|Y\setminus
F_{t-1}|}}\right)\\ &&\times
\left({\frac{v}{1-v}+\frac{1}{\psi_{|Y\setminus
F_{t-1}|}(q)}}\right)^{-1}\left({\frac{v}{1-v} +\frac{1}{\psi_{|F
_{t+1}\setminus Y|}(q)}}\right)^{-1}.
\end{eqnarray*}
Rewriting (\ref{star}) as $$ \left[\sum_{\cal F}\pi_{\cal F}{\cal
B}_{\cal F}(X, u)(1-u)\right]_{u=1}=\Phi_q(X;q,1), $$ we get
\begin{eqnarray*}
c(Y)&=&\left[\sum_{\cal F}\pi_{\cal F}{\cal B}_{\cal
F}(X;u)(1-p(Y)u)\right]_{u=p(-Y)}\\
&=&\Phi_q(-Y;q,1)\Phi_q(Y';q,1)\prod_{x_i\in Y, x_j\in
Y'}\frac{1-qx_i^{-1}x_j}{1-x_i^{-1}x_j}.
\end{eqnarray*}
Finally, the proof is completed by substituting  the values of
$c(Y)$ in (\ref{rational}).
\subsection{Some direct consequences  on $q$-series}
The following corollary of Theorem~1 will be usefull for the proof of
identities of Rogers-Ramanujan type.
\begin{theo} For $k\geq 1$,
\begin{eqnarray}
&&\hskip -1.5cm\sum_{l(\la)\leq
k}c_{(2\la)',k}(q)z^{|\la|}q^{n(2\la)}\left[{n\atop
2\la}\right]=(z;q^2)_n\sum_{r\geq 0}z^{kr}q^{(k+1)\left({2r \atop
2}\right)}\nonumber\\ &&\hskip 2cm \times \left[{n\atop
2r}\right]{1-zq^{4r-1}\over (zq^{2r-1})_{n+1}}.\label{eq:moi}\\
&&\hskip -1.5cm\sum_{l(\la)\leq
k}d_{\la',\,k}(q)z^{|\la|}q^{n(\la)}\left[{n\atop
\la}\right]=(z^2;q^2)_n\sum_{r\geq 0}z^{kr}q^{r+(k+1)\left({r\atop
2}\right)}\nonumber \\ &&\hskip 2cm\times \left[{n\atop
r}\right]{(1-zq^{-1})(1-z^2q^{2r-1})(1-zq^n)\over
(1-zq^{r-1})(1-zq^r)(z^2q^{r-1})_{n+1}}.\label{eq:autre}
\end{eqnarray}
\end{theo}
\begin{proof}
We know \cite[p. 213]{Ma} that if $x_i=z^{1/2}q^{i-1}$ ($1\leq
i\leq n$) then~:
\begin{equation}\label{value1}
P_{\la'}(X,q)=z^{|\la|/2}q^{n(\la)}\left[{n\atop \la}\right]. 
\end{equation}
Replacing $\la$ by $2\la$ and taking the conjugation in 
the left-hand side of (\ref{moi1}), we obtain  left-hand side of
(\ref{eq:moi}). On the other hand, for any $\xi\in \{\pm 1\}^n$ such
that the number of $\xi_i=-1$ is $r$, $0\leq r\leq n$,  we
have 
\begin{equation}\label{cal2}
\Phi_q(X^{\xi};0,0)= \Psi_q(X^{\xi};-1)\,\prod_i(1-x_i^{2\xi_i}),
\end{equation}
which is readily seen to equal 0 unless $\xi\in \{-1\}^r\times
\{1\}^{n-r}$. Now, in the latter case, 
 we have $\prod_ix_i^{k(1-\xi_i)/2}=z^{kr/2}q^{k\left({r
\atop 2}\right)}$, 
\begin{equation}\label{cal3}
\prod_{i=1}^n(1-x_i^{2\xi_i})=(-1)^rz^{-r}q^{-2\left({r \atop
2}\right)}(z;q^2)_{n},
\end{equation}
and \cite[p. 476]{St}~: 
\begin{equation}\label{cal4}
\Psi_q(X^{\xi};-1)=(-1)^rz^rq^{3\left({r \atop
2}\right)}\left[{n\atop r}\right]{1-zq^{2r-1}\over
(zq^{r-1})_{n+1}}. 
\end{equation}
Substituting these into 
the right side of 
(\ref{moi1}) with $r$ replaced by $2r$ we obtain the right side of
(\ref{eq:moi}). 

Similarly, in (\ref{autre}), replacing $x_i$ by $zq^{i-1}$ ($1\leq i\leq n$)
and invoking (\ref{value1}) we see that
the left side of (8) reduces to that of (\ref{eq:autre}). On the other
hand, since
$$
\Phi_q(X^{\xi};q,1)=\Phi_q(X^{\xi};0,0)\prod_{i=1}^n{1-qx_i^{\xi_i}\over
1-x_i^{\xi_i}},
$$
by (\ref{cal2}), this is equal to zero unless $\xi\in \{-1\}^r\times
\{1\}^{n-r}$ for some $r$, $0\leq r\leq n$. In the latter case, we
have
\begin{equation}
\prod_{i=1}^n{1-qx_i^{\xi_i}\over
1-x_i^{\xi_i}}=q^r{1-zq^{-1}\over 1-zq^{r-1}}{1-zq^n\over 1-zq^r},
\end{equation}
and invoking (\ref{cal2}), (\ref{cal3}) and (\ref{cal4}) with $z$
replaced by $z^2$,
\begin{equation}
\Phi_q(X^{\xi};0,0)  =q^{\left({r\atop 2}\right)}\left[{n\atop
r}\right](1-z^2q^{2r-1})
{(z^2;q^2)_n\over(z^2q^{r-1})_{n+1}}
\end{equation}
Plunging these into the right side of   (\ref{autre}) yields that of (\ref{eq:autre}).
\end{proof}\medskip

When $n\to +\infty$, Eqs. (\ref{eq:moi}) and  (\ref{eq:autre}) reduce
respectively to~:
\begin{equation}\label{lim1}
\sum_{l(\la)\leq k}\frac{z^{|\la|}q^{n(2\la)}}
{(q^2;q^2)_\la(q;q^2)_{\la_k}}
=(z;q^2)_\infty
\sum_{r\geq 0}{z^{kr}q^{(k+1)\bi{2r}{2}}\over 
(q)_{2r}(zq^{2r-1})_\infty}(1-zq^{4r-1}),
\end{equation}
\begin{eqnarray}
&&\sum_{l(\la)\leq k}
{z^{|\la|}q^{n(\la)}\over 
(q)_{\la_k}\,\prod_{i=1}^{k-1}(q^2;q^2)_{\left[(\la_i-\la_{i+1})/2\right]}}\label{lim2}\\ 
&=& (z^2;q^2)_\infty
\sum_{r\geq 0}z^{kr}q^{r+(k+1)\bi{r}{2}}
{1-zq^{-1}\over (q)_r(1-zq^{r-1})}
{1-z^2q^{2r-1}\over(1-zq^r)(z^2q^{r-1})_\infty}.\nonumber
\end{eqnarray}
Furthermore, setting  $z=q$ in  
(\ref{lim1}) and (\ref{lim2}) we obtain 
respectively  (11)
and \begin{equation}
\sum_{l(\la)\leq k}{q^{|\la|+n(\la)}\over (q)_{\la_k}\,
\prod_{i=1}^{k-1}(q^2;q^2)_{\left[(\la_i-\la_{i+1})/2\right]}}={1\over (q;q^2)_\infty}.
\end{equation}

\section{Elementary approach to multiple $q$-identities }
\subsection{Preliminaries}
Recall \cite[pp. 36-37]{An76} that the binomial formula has the
following $q$-analog~:
\begin{equation}
(z)_n=\sum_{m=0}^n\qbi{n}{m}{}(-1)^mz^mq^{m(m-1)/2}.\label{a1}
\end{equation}
Since the \emph{elementary symetric functions} $e_r(X)$ ($0\leq
r\leq n$) satisfy $$ (1+x_1z)(1+x_2z)\cdots
(1+x_nz)=\sum_{r=0}^ne_r(X)z^r, $$ it follows from (\ref{a1}) that
for integers $i\geq 0$ and $j\geq 1$~:
\begin{equation}\label{el}
e_r(q^i,\, q^{i+1},\, \ldots, \,q^{i+j-1})=q^{ir}e_r(1,\,
q,\,\ldots,\, q^{j-1})=q^{ir+\bi{r}{2}}\qbi{j}{r}{}.
\end{equation}
The following result can be derived from the Pieri's rule for
Hall-Littlewood polynomials~\cite[p. 215]{Ma}, but our proof is
elementary.
\begin{lem} For any partition $\mu$ such that $\mu_1\leq n$ there holds
\begin{equation}\label{qpieri}
q^{\bi{m}{2}+n(\mu)}\qbi{n}{m}{}\qbi{n}{\mu}{}=
\sum_{\la}q^{n(\la)}\qbi{n}{\la}{} \prod_{i\geq
1}\qbi{\la_i-\la_{i+1}}{\la_i-\mu_i}{},
\end{equation}
where the sum is over all partitions $\la$ such that
$\la/\mu$ is an  $m$-horizontal strip, i.e., $\mu\subseteq\la$,
$|\la/\mu|=m$ and there is at most one cell in each column of the
Ferrers diagram of $\la/\mu$.
\end{lem}
\begin{proof}
Let $l:=l(\mu)$ and $\mu_0=n$. Partition the set $\{1,\, 2,\,
\ldots, \, n\}$ into $l+1$ subsets~: $$ X_i=\{j\mid 1\leq j\leq
n\;\hbox{and}\;\mu_j'=i\}=\{j\mid \mu_{i+1}+1\leq j\leq \mu_i\},
\qquad 0\leq i\leq l.$$ Using (\ref{el}) to extract the
coefficients of $z^m$ in the following identity~: $$
(1+z)(1+zq)\cdots (1+zq^{n-1})=\prod_{i=0}^l\prod_{j\in
X_i}(1+zq^{j-1}), $$ we obtain
\begin{equation}\label{eq1}
q^{\bi{m}{2}}\qbi{n}{m}{}=\sum_{\bf
r}\prod_{i=0}^lq^{r_i\,\mu_{i+1}+\bi{r_i}{2}}\qbi{\mu_i-\mu_{i+1}}{r_i}{},
\end{equation}
where ${\bf r}=(r_0,\, r_1,\,\ldots, \, r_l)$ is a sequence of non
negative integers such that $\sum_ir_i=m$. For any such ${\bf r}$
define a partition $\la=(\la_1,\, \la_2,\, \ldots)$ by $$
\la_i=\mu_i+r_{i-1},\qquad 1\leq i\leq l+1. $$ Then $\la/\mu$ is a
$m$-horizontal strip. So (\ref{eq1}) can be written as
\begin{equation}\label{eq2}
q^{\bi{m}{2}}\qbi{n}{m}{}=\sum_{\la}\prod_{i=0}^lq^{(\la_{i+1}-\mu_{i+1})\mu
_{i+1}+\bi{\la_{i+1}-\mu_{i+1}}{2}}\qbi{\mu_i-\mu_{i+1}}{\mu_i-\la_{i+1}}{},
\end{equation}
where the sum is over all partitions $\la$ such that $\la/\mu$ is
a $m$-horizontal strip. Now, since
$$
(\la_{i+1}-\mu_{i+1})\mu_{i+1}+\bi{\la_{i+1}-\mu_{i+1}}{2}+\bi{\mu_{i+1}}{2}
=\bi{\la_{i+1}}{2},\qquad 0\leq i\leq l,
$$ 
and 
$\qbi{n}{\mu}{}
\prod_{i=0}^l\qbi{\mu_i-\mu_{i+1}}{\mu_i-\la_{i+1}}{}$ and
$\qbi{n}{\la}{} \prod_{i\geq
1}\qbi{\la_i-\la_{i+1}}{\la_i-\mu_i}{}$ are equal because they are both
equal to 
$$
{(q)_n\over (q)_{n-\la_1}
(q)_{\la_1-\mu_1}(q)_{\mu_1-\la_2}\cdots(q)_{\mu_l}},$$ 
 multiplying (\ref{eq2}) by
$q^{n(\mu)}\qbi{n}{\mu}{}$ yields (\ref{qpieri}).
\end{proof}

\begin{lem} There hold the following identities~:
\begin{eqnarray}
\sum_\la
z^{|\la|}q^{2n(\la)}\qbi{n}{\la}{}&=&\frac{1}{(z)_n}\label{hall},\\
\sum_\la
z^{|\la|}q^{n(\la)}\qbi{n}{\la}{}&=&\frac{(-z)_n}{(z^2)_n}\label{mac2},\\
\sum_\la (q,q^2)_\la\;
z^{|\la|}q^{n(2\la)}\qbi{n}{2\la}{}&=&\frac{(z;q^2)_n}{(z)_n}\label{mac3}.
\end{eqnarray}
\end{lem}
\begin{proof}
Identity (\ref{hall}) is due to Hall~\cite{H} and can be proved by
using the $q$-binomial identity~\cite{Ma2}. Stembridge~\cite{St}
proved (\ref{mac2}) using the $q$-binomial identity. Now, writing
$$\frac{(z^2;\, q^2)_n}{(z^2)_n}=(z)_n\,\frac{(-z)_n}{(z^2)_n}$$
 and applying successively
 (\ref{a1}), (\ref{mac2}) and
(\ref{qpieri}) we obtain
\begin{eqnarray*}
\frac{(z^2;\, q^2)_n}{(z^2)_n} &=&\sum_{\mu,
m}(-1)^mz^{m+|\mu|}q^{\bi{m}{2}+n(\mu)}
\qbi{n}{m}{}\qbi{n}{\mu}{}\\ &=& \sum_{\mu, m}(-1)^mz^{m+|\mu|}
\sum_{\la:\, \la/\mu=m-hs} \prod_{i\geq
1}\qbi{\la_i-\la_{i+1}}{\la_i-\mu_i}{}\,
q^{n(\la)}\qbi{n}{\la}{}\\
&=&\sum_{\la}z^{|\la|}q^{n(\la)}\qbi{n}{\la}{}\prod_{i\geq
1}\sum_{r_i\geq 0}(-1)^{r_i}\qbi{\la_i-\la_{i+1}}{r_i}{}.
\end{eqnarray*}
The identity (\ref{mac3}) follows then from $$
\sum_{j=0}^m(-1)^j\qbi{m}{j}{}=\left\lbrace\begin{array}{ll}
(q;\,q^2)_n&\quad \hbox{if}\quad m=2n,\\
 0&\quad \hbox{if $m$ is odd,}
\end{array}\right.
$$ which can be proved using the $q$-binomial formula~\cite[p.
36]{An76}.
\end{proof}

\begin{rem} When $n\to \infty$ the above identities reduce respectively to
the following~:
\begin{eqnarray}
\sum_\la
\frac{z^{|\la|}q^{2n(\la)}}{(q)_\la}&=&\frac{1}{(z)_\infty}\label{st1},\\
\sum_\la
\frac{z^{|\la|}q^{n(\la)}}{(q)_\la}&=&\frac{(-z)_\infty}{(z^2)_\infty}
\label{st2},\\ \sum_\la
\frac{z^{|\la|}q^{n(2\la)}}{(q^2;q^2)_\la}&=&\frac{1}{(zq;q^2)_\infty}\label
{st3}.
\end{eqnarray}
Also (\ref{st1}) and (\ref{st3}) are actually equivalent since the
later can be derived from (\ref{st1}) by substituting $q$ by $q^2$
and $z$ by $zq$.
\end{rem}

The following is the $q$-Gauss sum~\cite[p.10]{GR} due to Heine~:
\begin{equation}\label{Heine}
  {}_2\phi_1\left({a,b \atop x};q;{x\over
  ab}\right):= \sum_{n=0}^\infty
\frac{(a)_n(b)_n}{(q)_n(x)_n}\left(\frac{x}{ab}\right)^n={(x/a, \,
x/b;\, q)_\infty\over
  (x,x/ab;\, q)_\infty}.
\end{equation}
\begin{lem} We have
\begin{equation}
\sum_\la z^{|\la|}q^{n(2\la)} \frac{(a,\, b;
q^{-2})_{\la_1}}{(q^2;q^2)_{\la}}= \frac{(azq,\, bzq;
q^2)_\infty}{(zq,\, abzq; q^2)_\infty}\label{d}.
\end{equation}
\end{lem}
\begin{proof}
 Substituting $q^2$ by $q$ and $z$ by
$zq$, the identity is equivalent to
\begin{equation}\label{gauss}
\sum_\la z^{|\la|}q^{2n(\la)} \frac{(a,\, b;
q^{-1})_{\la_1}}{(q)_{\la}}= \frac{(az,\, bz; q)_\infty}{(z,\,
abz; q)_\infty}.
\end{equation}
Now, writing $k=\la_1$ and $\mu=(\la_2,\la_3,\cdots)$, and using
(\ref{hall}) we get
\begin{eqnarray*}
\sum_\la z^{|\la|}q^{2n(\la)} \frac{(a,\, b;
q^{-1})_{\la_1}}{(q)_{\la}}&=&\sum_{k\geq
0}z^kq^{k(k-1)}{(a,b;q^{-1})_k\over (q)_k}\sum_\mu
z^{|\mu|}q^{2n(\mu)}\qbi{k}{\mu}{}\\ &=&\sum_{k\geq
0}(abz)^k{(a^{-1},b^{-1}; q)_k\over (q)_k(z)_k}.
\end{eqnarray*}
Identity (\ref{gauss}) follows then from (\ref{Heine}).
\end{proof}

\begin{rem}
Formula~(\ref{gauss}) was derived in \cite{St} from a more general
formula of Hall-Littlewood polynomials.
\end{rem}

\subsection{Elementary proof of Theorem~4}
We shall only prove (\ref{eq:moi}) when $n$ is even 
and leave the case when $n$ is odd  and 
 (\ref{eq:autre}) to the
interested reader because their proofs are very similar.
  Consider the
generating function of the left-hand side of (\ref{eq:moi}) with $n=2r$~:
\begin{eqnarray}
\varphi(u)
&=&\sum_{k\geq 0}u^k\sum_{l(\la)\leq
k}\frac{(q;q^2)_\la}{(q;q^2)_{\la_k}}z^{|\la|}q^{n(2\la)}\left[{2r\atop
2\la}\right]\nonumber\\ 
&=&\sum_{\la}
u^{l(\la)}z^{|\la|}q^{n(2\la)}(q;q^2)_\la
\left[{2r\atop 2\la}\right]\sum_{k\geq 0} {u^k\over (q;q^2)_{\la_{k+l(\la)}}}\nonumber\\ 
&=&\sum_{\la}u^{l(\la)}z^{|\la|}q^{n(2\la)}(q;q^2)_\la
\left[{2r\atop
2\la}\right]\left(\frac{u}{1-u}+\frac{1}{(q;q^2)_{\la_{l(\la)}}}\right)
\label{ref1}.
\end{eqnarray}
Now, each  
partition $\la$ with parts bounded by $r$ can be encoded by
 a pair of sequences
 $\nu=(\nu_0, \nu_1, \cdots, \nu_l)$ and ${\mathbf m}=(m_0,\cdots,m_{l})$
such that $\la=(\nu_0^{m_0},\ldots, \nu_l^{m_l})$, where
$r=\nu_0>\nu_1> \cdots
>\nu_l>0$ and $\nu_i$ has multiplicity $m_i\geq 1$ for $1\leq
i\leq l$ and $\nu_0=r$ has   multiplicity $m_0\geq 0$. 
Using the notation~: 
$$
<\alpha>=\frac{\alpha}{1-\alpha},\quad u_i=z^iq^{i(2i-1)}\quad\hbox{for}\quad
 i\geq 0,
$$ 
we can then rewrite (\ref{ref1})
as follows~:
\begin{eqnarray}
\varphi(u)&=&\sum_{\nu}(q;q^2)_\nu\left[{2r\atop
2\nu}\right]\left(<u>+\frac{1}{(q;q^2)_{\nu_l }}\right)\nonumber\\
&&\hskip 10pt \times\sum_{{\mathbf m}}\left((u_ru)^{m_0}+{\chi(m_0=0)\over
(q;q^2)_{r-\nu_1}}\right)\prod_{i=1}^{l}(u_{\nu_i}u)^{m_i
}\nonumber\\ 
&=&\sum_{\nu}{(q)_{2r}\over (q^2;q^2)_\nu}B_\nu,\label{x}
\end{eqnarray}
where the sum is over all strict partitions $\nu=(\nu_0,\nu_1,\ldots, \nu_l)$ and
$$
B_\nu=\left(<u>+\frac{1}{(q;q^2)_{\nu_l
}}\right)\left( <u_ru>+\frac{1}{(q;q^2)_{r-\nu_1}}\right) 
\prod_{i=1}^{l}<u_{\nu_i}u>.
$$
So $\varphi(u)$ is  a rational fraction with simple poles at 
$u_p^{-1}$ for $0\leq p\leq r$. Let $b_p(z,r)$ be the corresponding residue of 
 $\varphi(u)$ at $u_p^{-1}$ for $0\leq p\leq r$.
Then, it follows from (\ref{x}) that
\begin{equation}\label{residue}
b_p(z,r)=\sum_\nu{(q)_{2r}\over
(q^2;q^2)_\nu}\left[B_\nu(1-u_pu)\right]_{u=u_p^{-1}}.
\end{equation}
We shall first consider the cases where $p=0$ or $r$.
Using
 (\ref{ref1}) and (\ref{mac3}) we have
\begin{equation}\label{b0}
b_0(z,r)=\left[\varphi(u)(1-u)\right]_{u=1}={(z;q^2)_{2r}\over (z)_{2r}}.
\end{equation}
Now, by (\ref{x}) and(\ref{residue}) we have
\begin{equation}\label{b0bis}
b_0(z,r)=\sum_{\nu}{(q)_{2r}\over (q^2;q^2)_\nu}
\left(<u_r>+\frac{1}{(q;q^2)_{r-\nu_1}}\right)
\prod_{i=1}^{l}<u_{\nu_i}>,
\end{equation}
and
$$
b_r(z,r)=\sum_{\nu}{(q)_{2r}\over (q^2;q^2)_\nu}\left(<1/u_r>+\frac{1}{(q;q^2)_{\nu_l
}}\right)
\prod_{i=1}^{l}<u_{\nu_i}/u_r>,
$$
which, by setting $\mu_i=r-\nu_{l+1-i}$ for $1\leq
i\leq l$ and $\mu_0=r$,   can be written as 
\begin{equation}\label{br}
b_r(z,r)=\sum_{\mu}{(q)_{2r}\over
(q^2;q^2)_\mu}\left(<1/u_r>+\frac{1}{(q;q^2)_{r-\mu_1
}}\right)\prod_{i=1}^{l}<u_{r-\mu_i}/u_r>.
\end{equation}
Comparing (\ref{br}) with (\ref{b0bis}) we see that
$b_r(z,r)$ is equal to $b_0(z,r)$ with $z$ replaced by $z^{-1}q^{-2(2r-1)}$. 
Il follows from
(\ref{b0}) that
\begin{equation}\label{ref3}
b_r(z,r)=b_0(z^{-1}q^{-2(2r-1)},r)=(z;q^2)_{2r}q^{r(2r-1)}\,\frac{1-zq^{4r-1}}
{(zq^{2r-1})_{2r+1}}. 
\end{equation}
Consider now the case where  $0<p<r$. Clearly,  
for each partition $\nu$, 
the corresponding
 summand in (\ref{residue}) is not zero only if 
$\nu_j=p$ for some $j$, $0\leq j\leq r$. 
Furthermore, each such partition $\nu$ can be splitted into two
strict partitions $\rho=(\rho_0, \rho_1,\ldots, \rho_{j-1})$ and 
$\sigma=(\sigma_0,
 \ldots, \sigma_{l-j})$ such that
$\rho_i=\nu_i-p$ for $0\leq i\leq j-1$ and $ \sigma_s=\nu_{j+s}$ for 
$0\leq s\leq l-j$. So we can write (\ref{residue}) as follows~:
\begin{eqnarray*}
b_p(z, r)&=&\left[{2r\atop 2p}\right]
\sum_{\rho}\frac{(q)_{2r-2p}}{(q^2;q^2
)_\rho}F_\rho(p)\times\sum_\sigma\frac{(q)_{2p}}{(q^2;q^2)_\sigma}G_\sigma(p)
\end{eqnarray*}
where  for $\rho=(\rho_0, \rho_1,\ldots, \rho_l)$ with $\rho_0=r-p$,
$$
F_\rho(p)=\left(<u_r/u_p>+{1\over
(q;q^2)_{r-p-\rho_1}}\right)\prod_{i=1}^{l(\rho)}<u_{\rho_i+p}/u_p>,
$$
and for $\sigma=(\sigma_0, \ldots, \sigma_l)$ with $\sigma_0=p$,
$$
G_\sigma(p)=\left(<1/u_p> +\frac{1}{(q;q^2)_{\sigma_{l} }}\right)
\prod_{i=1}^{l(\sigma)}<u_{\sigma_i}/u_p>.
$$
Comparing with (\ref{b0bis}) and (\ref{br}) and 
using (\ref{b0}) and (\ref{ref3}) we obtain
\begin{eqnarray*}
b_p(z,r)&=&\left[{2r\atop 2p}\right]b_0(zq^{4p},r-p)\,b_p(z,p)\\
 &=&\left[{2r\atop
2p}\right](z;q^2)_{2r}q^{\left({2r\atop
2p}\right)}\frac{1-zq^{4p-1}}{(zq^{2p-1})_{2r+1}}.
\end{eqnarray*}
Finally,  extracting the
coefficients of $u^k$ in the equation
$$
\varphi(u)=\sum_{p=0}^r\frac{b_p(z,r)}{1-u_pu},
$$ 
and using  the values for $b_p(z,r)$ we obtain(\ref{eq:moi}).

\subsection{Proof of Theorem 2}
Consider the generating function of the left-hand side of
(\ref{moi5}) :
\begin{eqnarray}
\varphi_{ab}(u)&:=&\sum_{k\geq 0}u^k\sum_{l(\la)\leq
k}z^{|\la|}q^{n(2\la)}\frac{(a,\,
b;q^{-2})_{\la_1}}{(q^2;q^2)_\la(q;q^2)_{\la_k}}\nonumber\\
&=&\sum_\la\sum_{k\geq
0}u^{k+l(\la)}z^{|\la|}q^{n(2\la)}\frac{(a,\,
b;q^{-2})_{\la_1}}{(q^2;q^2)_\la(q;q^2)_{\la_{l(\la)+k}}}\nonumber\\
&=&\sum_\la u^{l(\la)}z^{|\la|}q^{n(2\la)}\frac{(a,\,
b;q^{-2})_{\la_1}}{(q^2;q^2)_\la}\left(\frac{u}{1-u}+\frac{1}{(q;q^2)_{\la_{
l(\la)}}}\right)\label{g'}.
\end{eqnarray}
As in the proof of Theorem~4, we encode the partition $\la$ in
the previous sum. Let $\nu_1, \cdots, \nu_l, \nu_{l+1}=0$ denote
the distinct parts of $\la$, so that $\nu_1> \cdots
>\nu_l>\nu_{l+1}=0$ and $\nu_i$ has multiplicity $m_i$ for $1\leq
i\leq l$. Then we have
\begin{eqnarray}
\varphi_{ab}(u)&=&\sum_{\nu,\,{\mathbf m}}\frac{(a,\,
b;q^{-2})_{\nu_1}}{(q^2;q^2)_\nu}\left(\frac{u}{1-u}+\frac{1}{(q;q^2)_{\nu_l
}}\right)
\prod_{i=1}^{l}(u_{\nu_i}u)^{m_i}\nonumber\\
&=&\sum_{\nu}\frac{(a,\,
b;q^{-2})_{\nu_1}}{(q^2;q^2)_\nu}\left(<u>+\frac{1}{(q;q^2)_{\nu_l}}\right)
\prod_{i=1}^l<u_{\nu_i}u>\label{h}.
\end{eqnarray}
Each of the terms in this sum, as a rational function of $u$, has
a finite set of simple poles, which may occur at the points
$u_r^{-1}$ for $r\geq 0$. Therefore, each term is a
linear combination of partial fractions. Moreover, the sum of
their expansions converges coefficientwise. So $\varphi_{ab}$ has
an expansion 
$$\varphi_{ab}(u)=\sum_{r\geq
0}\frac{c_r}{1-uz^rq^{r(2r-1)}},
$$ where $c_r$ denotes the formal
sum of partial fraction coefficients contributed by the terms of
(\ref{h}). It remains to compute these residues $c_r$ $(r\geq 0)$.
By using (\ref{d}) and (\ref{g'}), we get immediately
$$c_0=\varphi_{ab}(u)(1-u)|_{u=1}= \frac{(azq,\, bzq;
q^2)_\infty}{(zq,\, abzq; q^2)_\infty}.$$ In view of (\ref{h}),
this yields the identity
\begin{equation}\label{i}
\sum_{\nu}\frac{(a,\,
b;q^{-2})_{\nu_1}}{(q^2;q^2)_\nu}\prod_{i=1}^{l}
<u_{\nu_i}>=\frac{(azq,\, bzq;
q^2)_\infty}{(zq,\, abzq; q^2)_\infty}.
\end{equation}
To compute the residues $c_r$ for $r>0$, the contribution in
(\ref{h}) is given only by the partitions $\nu$ for which $\exists
j\, |\,\nu_j=r$. For each such partition, we define as before
$\rho_i:=\nu_i-r$ for $1\leq i<j$ and $\sigma_i:=\nu_{i+j}$ for
$0\leq i\leq l-j$. We get two partitions $\rho$ and $\sigma$ with
$\sigma$ bounded by $r$. Using (\ref{h}), we obtain
\begin{eqnarray*}
c_r&=&\left[\varphi_{ab}(u)(1-u_ru)\right]_{u=u_r^{-1}}\\
&=&\sum_{\rho,\,\sigma}\frac{(a,b;q^{-2})_{\rho_1+r}}{(q^2;q^2)_\rho(q^2;q^2
)_\sigma}
\prod_{i=1}^{j-1}<u_{r+\rho_i}/u_r>\\
&&\hskip 1 cm
\times\left(<1/u_r>+\frac{1}{(q;q^2)_{\sigma_l}}\right)
\prod_{i=1}^{l-j}<u_{\sigma_i}/u_r>.
\end{eqnarray*}
To eliminate the $\sigma$-dependence of this series, we apply
(\ref{ref3}), and this leads to
\begin{eqnarray*}
c_r&=&\sum_{\rho}\frac{(aq^{-2r},bq^{-2r};q^{-2})_{\rho_1}(a,b;q^{-2})_r}
{(q^2;q^2)_\rho(q)_{2r}}\\ 
&&\hskip 2cm
\times(z;q^2)_{2r}q^{\left({2r\atop
2}\right)}\frac{1-zq^{4r-1}}{(zq^{2r-1})_{2r+1}}\prod_{i=1}^{j-1}<{u_{r+\rho_i}\over u_r}>\\
&=&\frac{(a,b;q^{-2})_r(z;q^2)_{2r}}{(q)_{2r}}q^{\left({2r\atop
2}\right)}\frac{1-zq^{4r-1}}{(zq^{2r-1})_{2r+1}}\frac{(azq^{-2r+1+4r},bzq^{
-2r+1+4r};q^2)_\infty} {(zq^{4r+1},abzq;q^2)_\infty},
\end{eqnarray*}
where the last equality follows from (\ref{i}) with $a$, $b$, $z$
replaced respectively by $aq^{-2r}$, $bq^{-2r}$, $zq^{4r}$. After
simplification, one gets 
$$
c_r=q^{\left({2r\atop
2}\right)}\frac{(z;q^2)_\infty}{(zq^{2r-1})_\infty}\frac{(a,b;q^{-2})_r
(azq^{2r+1},bzq^{2r+1};q^2)_\infty}{(q)_{2r}(abzq;q^2)_\infty}(1-zq^{4r-1}
), $$ which completes the proof.

\section{Comparison with Andrews-Bailey's method}
A popular method  to prove identities of Rogers-Ramanujan type
is based on Bailey's lemma (see \cite{An86,Wa}).
In~\cite{An84} Andrews noticed that by  appling
iteratively Bailey's lemma to the corresponding Bailey pair in the
simple sum case one can obtain 
multianalog identities of 
Rogers-Ramanujan type almost straightforwardly.
In this section we shall 
briefly
compare our multisum analogs with those obtained through Andrews-Bailey's
approach. Recall that a pair
 of sequences $(\alpha_n)_{n\geq 0}$ and $(\beta_n)_{n\geq 0}$ is a \emph{Bailey pair}
 if they are
related by the following~\cite[p. 25-26]{An86}~:
\begin{equation}\label{pair}
\beta_n=\sum_{r=0}^n{\alpha_r\over (q)_{n-r}(aq)_{n+r}}.
\end{equation}
If $(\alpha_n,\,\beta_n)$ is a Bailey
pair  and $(\alpha'_n,\,\beta'_n)$ is one of the
following pairs~:
\begin{eqnarray*}
&(i)&\quad \alpha'_n=a^nq^{n^2}\alpha_n,\qquad  \beta'_n=\sum_{k\geq 0}{a^kq^{k^2}\over
(q)_{n-k}}\beta_k;\label{Ba1}\\ 
&(ii)&\quad \alpha'_n={(-q^{1/2})_n\over
(-aq^{1/2})_n}a^nq^{n^2/2}\alpha_n,\qquad \beta'_n=\sum_{k\geq 0}{(-q^{1/2})_ka^kq^{k^2/2}\over
(q)_{n-k}(-aq^{1/2})_k}\beta_k;\label{Ba2}\\
&(iii)&\quad \alpha'_n=a^{n/2}q^{n^2/2}\alpha_n,\qquad \beta'_n=\sum_{k\geq
0}{(-(aq)^{1/2})_ka^{k/2}q^{k^2/2}\over
(q)_{n-k}(-(aq)^{1/2})_n}\beta_k\label{Baa2};
\end{eqnarray*}
 then  Bailey's lemma~\cite[p. 25-26]{An86} states that
$(\alpha'_n,\,\beta'_n)$ is  also a Bailey pair. What we need
here is actually the limit case of (\ref{pair}). 
In (\ref{pair}) substituting
$(\alpha_n,\, \beta_n)$ by one of the above $(\alpha'_n,\,\beta'_n)$'s
and letting  $n\to\infty$,  we obtain respectively
\begin{eqnarray}
&&\hskip -1cm\sum_{n\geq 0}a^nq^{n^2}\beta_n={1\over
(aq)_\infty}\sum_{r\geq 0}a^rq^{r^2}\alpha_r\label{Ba3},\\
&&\hskip -1cm\sum_{n\geq
0}a^nq^{n^2/2}(-q^{1/2})_n\beta_n={(-aq^{1/2})_\infty\over(aq)_\infty}
\sum_{r\geq
0}{(-q^{1/2})_r\over
(-aq^{1/2})_r}a^rq^{r^2/2}\alpha_r\label{Ba4},\\ &&\hskip
-1cm\sum_{n\geq
0}a^{n/2}q^{n^2/2}(-(aq)^{1/2})_n\beta_n
={(-(aq)^{1/2})_\infty\over(aq)_\infty}\sum_{r\geq
0}a^{r/2}q^{r^2/2}\alpha_r\label{Ba5}.
\end{eqnarray}
Now, if we iterate the above process
 $k$ times to a same Bailey pair \cite[p.30]{An86}, then
(\ref{Ba3}), (\ref{Ba4}) and (\ref{Ba5}) lead  respectively to
\begin{eqnarray}
&&\hskip -1cm{1\over (aq)_\infty}\sum_{r\geq
0}a^{kr}q^{kr^2}\alpha_r=\sum_{l(\la)\leq
k}a^{|\la|}q^{n_2(\la)}{(q)_{\la_k}\over
(q)_\la}\beta_{\la_k},\label{an1}\\ &&\hskip
-1cm{(-q^{1/2})_\infty\over (aq)_\infty}\sum_{r\geq
0}\left({(-q^{1/2})_r\over
(-aq^{1/2})_r}\right)^ka^{kr}q^{kr^2/2}\alpha_r\label{an2}\\
&&\hskip 1cm=\sum_{l(\la)\leq
k}a^{|\la|}q^{n_2(\la)/2}{(-q^{1/2})_{\la_1}\cdots(-q^{1/2})_{\la_k}\over
(-aq^{1/2})_{\la_1}\cdots(-aq^{1/2})_{\la_{k-1}}}{(q)_{\la_k}\over
(q)_\la}\beta_{\la_k},\nonumber\\&&\hskip -1cm
{(-(aq)^{1/2})_\infty\over (aq)_\infty}\sum_{r\geq
0}a^{kr/2}q^{kr^2/2}\alpha_r\label{an3}\\&&\hskip
1cm=\sum_{l(\la)\leq
k}a^{|\la|/2}q^{n_2(\la)/2}{(-(aq)^{1/2})_{\la_k}(q)_{\la_k}\over
(q)_\la}\beta_{\la_k}.\nonumber
\end{eqnarray}

Slater~\cite{Sl1, Sl2} noticed that (\ref{moi3}) and (\ref{moi9}) follow
from (\ref{Ba3}) and (\ref{Ba4}) by choosing the  pair :
\begin{equation}\label{pair1}
\alpha_n=q^{n^2}(q^{n/2}+q^{-n/2}),\qquad \beta_n={1\over
(q^{1/2}, q;\, q)_n},
\end{equation}
with $a=1$ and $q$ replaced by $q^2$, and 
(\ref{moi4}) and (\ref{moi10}) follow from 
(\ref{Ba3}) and (\ref{Ba5}) by choosing  the pair :
\begin{equation}\label{pair2}
\alpha_n=q^{n^2+n/2}{1+q^{n+{1/2}}\over 1+q^{1/2}}, \qquad
\beta_n={1\over (q^{3/2}, q;\, q)_n}
\end{equation}
with $a=q$ and $q$ replaced by $q^2$.

Now, if we choose the Bailey pair (\ref{pair1}) in (\ref{an1})
and (\ref{an2}) with $a=1$ and $q$ replaced by $q^2$, we obtain respectively
(11) and  
\begin{equation}\label{new1}
\sum_{l(\la)\leq k}
\frac{q^{n_2(\lambda)}(-q;q^2)_{\lambda_k}}
{(q;q^2)_{\lambda_k}(q^2;q^2)_{\la}}=
\frac{(q^{2k+4},\,
-q^{k+1},-q^{k+3};\;q^{2k+4})_\infty}{(q)_\infty\,(-q^2;q^2)_\infty},
\end{equation}
which is different from
(\ref{moi7}).
In the same way, if we choose the
pair (\ref{pair2}) in 
(\ref{an1}) and (\ref{an3}) with $a=q$ and $q$ replaced by $q^2$ 
then we obtain
\begin{eqnarray}\label{new2}
\sum_{l(\la)\leq k}
\frac{q^{2|\la|+2n_2(\lambda)}}{(q;q^2)_{\lambda_k}(q^2;q^2)_{\lambda}}&=&
\frac{(q^{4k+4},\,
-q^{4k+3},-q;\;q^{4k+4})_\infty}{(q^2;q^2)_\infty},\\
\sum_{l(\la)\leq k}
\frac{q^{|\la|+n_2(\lambda)}(-q^2;q^2)_{\lambda_k}}{(q;q^2)_{\lambda_k}(q^2;
q^2)_{\lambda}}&=&
\frac{(q^{2k+4},\,
-q^{2k+3},-q;\;q^{2k+4})_\infty}{(q)_\infty\,(-q;q^2)_\infty}.\label{new3}
\end{eqnarray}
which are different from
(\ref{moi6}) and (\ref{moi8}), respectively.

So, only equation (11) coincides with that directly obtained by Andrews-Bailey's method. It seems that some new techniques may be necessary
to demonstrate  all our six multisum identities of
Rogers-Ramanujan type  through the classical Andrews-Bailey's 
method. Recently, Bressoud, Ismail and Stanton~\cite{BIS} have proved all the
sixteen identities in Stembridge's paper~\cite{St} by means of 
\emph{change of base in Bailey pairs}. It 
would be interesting to see whether their method can be applied to
our identities.



\small
\begin{thebibliography}{99}
\bibitem{ACE} \textsc{Veigneau} {S.}, \emph{ACE, an Algebraic Environment
for the Computer algebra system MAPLE}, \texttt{http:
file://phalanstere.univ-mlv.fr/~ace}, 1998.
\bibitem{An76} \textsc{Andrews} (G.E.), \emph{The theory of partitions},
Encyclopedia of mathematics and its applications,
 Vol. {\bf 2}, Addison-Wesley, Reading, Massachusetts, 1976.
%
\bibitem{An84} \textsc{Andrews} (G.E.),
Multiple series Rogers-Ramanujan type identities, \emph{Pacific J.
Math.}, Vol. {\bf 114}, No. 2, 267-283, 1984.
\bibitem{An86} \textsc{Andrews} (G.E.),
$q$-Series~: Their Development and Application in Analysis,
Combinatorics, Physics, and Computer Algebra, \emph{CBMS Regional
Conference Series},  Vol. {\bf 66}, Amer. Math. Soc., Providence,
1986.

\bibitem{BIS} \textsc{Bressoud} (D.), \textsc{Ismail} (M.), and
\textsc{Stanton} (D.), Change of Base in Bailey Pairs, \emph{The
Ramanujan J.}, {\bf 4}, 435-453, 2000.
\bibitem{GR} \textsc{Gasper} (G.) and \textsc{Rahman} (M.),
\emph{Basic Hypergeometric Series}, Encyclopedia of mathematics
and its applications, 35, Cambridge University Press, Cambridge,
1990.
\bibitem{H} \textsc{Hall} (P.), A partition formula connected
with Abelian groups, \emph{Comment. Math. Helv.} {\bf 11},
126-129, 1938.
\bibitem{JZ} \textsc{Jouhet} (F.) and \textsc{Zeng} (J.), Some
 new identities for Schur functions, to appear in \emph{Adv. Appl. Math.},
2001.
\bibitem{Ma2} \textsc{Macdonald} (I.G.),
An elementary proof of a $q$-binomial identity, \emph{$q$-series
and partitions} (Minneapolis, MN, 1988), 73--75, IMA Vol. Math.
Appl., 18, Springer, New York, 1989.
\bibitem{Ma} \textsc{Macdonald} (I.G.),
\emph{Symmetric functions and Hall polynomials}, Clarendon Press,
second edition, Oxford, 1995.
\bibitem{Sl1} \textsc{Slater} (L. J.), A new proof of Rogers's
transformations of infinite series, \emph{Proc. London Math.
Soc.}, {\bf 53} (2), 460-475, 1951.
\bibitem{Sl2} \textsc{Slater} (L. J.), Further identities of
the  Rogers-Ramanujan Type, \emph{Proc. London  Math. Soc.}, {\bf
54} (2), 147-167 (1951-52).
\bibitem{St} \textsc{Stembridge} (J. R.), Hall-Littlewood
functions, plane partitions, and the Rogers-Ramanujan identities,
\emph{Trans. Amer. Math. Soc.}, {\bf 319}, no.2, 469-498, 1990.
\bibitem{Wa} \textsc{Warnaar} (S. O.), 50 years of Bailey's lemma, 
\emph{Algebraic combinatorics and applications},  Springer-Verlag,
2001.
\end{thebibliography}
\end{document}